\documentclass[graybox]{svmult}

\usepackage{amsmath,amssymb}
\usepackage{mathptmx}       
\usepackage{helvet}         
\usepackage{courier}        
\usepackage{type1cm}        
\usepackage{makeidx}         
\usepackage{multicol}        
\usepackage[bottom]{footmisc}

\allowdisplaybreaks

\newtheorem{dfn}{Definition}
\newtheorem{thm}[dfn]{Theorem}

\newtheorem{cor}[dfn]{Corollary}
\newtheorem{lem}[dfn]{Lemma}
\newtheorem{rem}[dfn]{Remark}
\newtheorem{pro}[dfn]{Proposition}

\newcommand{\R}{{\mathbb{R}}}

\renewcommand{\a}{\alpha}
\renewcommand{\b}{\beta}
\renewcommand{\l}{\lambda}

\makeindex

\begin{document}
	
\title*{New inequalities for $\eta$-quasiconvex functions\thanks{This 
is a preprint of a paper whose final and definite form is  
accepted 19-Sept-2018 as a book chapter at Springer New York, 
on the topic of ``Frontiers in Functional Equations and Analytic Inequalities'',
Edited by G. Anastassiou and J. Rassias.}}

\author{Eze R. Nwaeze \and Delfim F. M. Torres}

\institute{Eze R. Nwaeze 
\at Department of Mathematics, Tuskegee University,
Tuskegee, AL 36088, USA\\ 
\email{enwaeze@tuskegee.edu}
\and Delfim F. M. Torres (corresponding author)
\at CIDMA, Department of Mathematics, University of Aveiro,
3810-193 Aveiro, Portugal\\ 
\email{delfim@ua.pt}}

\maketitle

\abstract{The class of $\eta$-quasiconvex functions was introduced in 2016. 
Here we establish novel inequalities of Ostrowski type for functions 
whose second derivative, in absolute value raised to the power $q\geq 1$, 
is $\eta$-quasiconvex. Several interesting inequalities are deduced as special cases. 
Furthermore, we apply our results to the arithmetic, geometric, Harmonic, 
logarithmic, generalized log and identric means, 
getting new relations amongst them.}

\bigskip


\medskip



\section{Introduction}
\label{intro}

A function $G:I\rightarrow\R$ is said to be convex on the interval 
$I\subset\R$ if 
$$
G(xu+(1-x)v)\leq xG(u)+(1-x)G(v)
$$
holds for all $u,v\in I$ and $x\in[0, 1]$.
Many interesting inequalities have been established 
for convex functions. Worthy of mention is the following 
result proved in 2011 by Sarikaya and Aktan \cite{SA}.

\begin{thm}[See \cite{SA}]
\label{thm:01}
Let $I\subset\R$ be an open interval, $\a,\b \in \R$ 
with $\a<\b$, $\l\in[0, 1]$, and $G:I\rightarrow\R$ 
be a twice differentiable mapping such 
that $G''$ is integrable. If $|G''|$ 
is a convex function on $[\a,\b]$, then 
\begin{align*}
&\Bigg|(\l-1)G\left(\frac{\a+\b}{2}\right)
-\l\frac{G(\a)+G(\b)}{2}+\frac{1}{\b-\a}\int_\a^\b G(x)\,dx\Bigg|\\
&\leq
\begin{cases}
\frac{(\b-\a)^2}{12}\Bigg[\Big(\l^4+(1+\l)(1-\l)^3+\frac{5\l-3}{4}\Big) |G''(\a)|\\
\qquad +\Big(\l^4+(2-\l)\l^3+\frac{1-3\l}{4}\Big)|G''(\b)|\Bigg],
\quad {\rm if}~~~0\leq\l \leq \frac{1}{2},\\
\frac{(\b-\a)^2(3\l-1)}{48}\Big[|G''(\a)|+|G''(\b)|\Big], 
\quad {\rm if}~~\frac{1}{2}\leq \l \leq 1.
\end{cases}
\end{align*}
\end{thm}

In 2015, Liu obtained a related inequality 
for $s$-convex functions \cite{Liu}. 
The notion of $s$-convexity was introduced 
in 1994 by Hudzik and Maligranda \cite{Hud}.
Let us recall it here.

\begin{dfn}[See \cite{Hud}]
\label{def:s:conv}
A function $G:[0,\infty)\rightarrow\R$ is said to be $s$-convex if 
$$
G(ux+(1-x)v)\leq x^sG(u)+(1-x)^sG(v)
$$
holds for all $u,v\in I,$ $x\in[0, 1]$ and for some fixed $s\in (0,1]$.
\end{dfn}

Evidently, the notion of $s$-convexity given in Definition~\ref{def:s:conv}
generalizes the classical concept of convexity. 
For this class of functions, Liu \cite{Liu}, 
among other things, established the following result:

\begin{thm}[See \cite{Liu}]
\label{thm:Liu}
Let $I\subset[0,\infty)$, $G:I\rightarrow\R$ be a twice differentiable 
function on $I^{\circ}$ such that  $G''\in L_1[\a,\b]$, where 
$\a,\b\in I$ with $\a<\b$. If $|G''|^q$ is $s$-convex on $[\a, \b]$ 
for some fixed $s\in(0, 1]$ and $q\geq1$, then 
\begin{align*}
\Bigg|&\frac{1}{\b-\a}\int_\a^\b G(x)\,dx
-(1-\l)G\Big(\frac{\a+\b}{2}\Big)-\l\frac{G(\a)+G(\b)}{2}\Bigg|\\
&\leq \frac{(\b-\a)^2}{16}\Big(\frac{8\l^3-3\l+1}{3}\Big)^{1-\frac{1}{q}}\Bigg\{\Bigg[
\frac{2(2\l)^{s+3}-2(s+3)\l+s+2}{(s+2)(s+3)}\Bigg|G''\Bigg(\frac{\a+\b}{2}\Bigg)\Bigg|^q\\
&\quad +\frac{4(1-2\l)^{s+2}[(s+1)\l+1]+2(s+3)\l-2}{(s+1)(s+2)(s+3)}|G''(\a)|^q\Bigg]^{\frac{1}{q}}\\
&\quad + \Bigg[ \frac{2(2\l)^{s+3}-2(s+3)\l+s+2}{(s+2)(s+3)}\Bigg|G''\Bigg(\frac{\a+\b}{2}\Bigg)\Bigg|^q\\
&\quad +\frac{4(1-2\l)^{s+2}[(s+1)\l+1]+2(s+3)\l-2}{(s+1)(s+2)(s+3)}|G''(\b)|^q\Bigg]^{\frac{1}{q}} \Bigg\}
\end{align*}
for $0\leq \l\leq \frac{1}{2}$ and
\begin{align*}
\Bigg|&\frac{1}{\b-\a}\int_\a^\b G(x)\,dx-(1-\l)
G\Big(\frac{\a+\b}{2}\Big)-\l\frac{G(\a)+G(\b)}{2}\Bigg|\\
&\leq \frac{(\b-\a)^2}{16}\Big(\l-\frac{1}{3}\Big)^{1-\frac{1}{q}}\Bigg\{\Bigg[
\frac{2(s+3)\l-s-2}{(s+2)(s+3)}\Bigg|G''\Bigg(\frac{\a+\b}{2}\Bigg)\Bigg|^q\\
&\quad +\frac{2(s+3)\l-2}{(s+1)(s+2)(s+3)}|G''(\a)|^q\Bigg]^{\frac{1}{q}}+ \Bigg[
\frac{2(s+3)\l-s-2}{(s+2)(s+3)}\Bigg|G''\Bigg(\frac{\a+\b}{2}\Bigg)\Bigg|^q\\
&\quad +\frac{2(s+3)\l-2}{(s+1)(s+2)(s+3)}|G''(\b)|^q\Bigg]^{\frac{1}{q}}\Bigg\}
\end{align*}
for $\frac{1}{2}\leq \l\leq 1$.
\end{thm}

In 2016, Eshaghi Gordji et al. \cite{GDS} proposed  a larger class 
of functions called $\eta$-quasiconvex. 

\begin{dfn}[See \cite{GDS}]
\label{eta}
A function $G:I\subset\mathbb{R}\rightarrow\mathbb{R}$ is said to be 
an $\eta$-quasiconvex function with respect to 
$\eta:\mathbb{R}\times\mathbb{R}\rightarrow\mathbb{R}$, if 
$$
G(xu+(1-x)v)\leq \max\left\{G(v), G(v)+\eta(G(u), G(v))\right\}
$$
for all $u, v\in I$ and $x\in [0, 1]$.
\end{dfn}

An $\eta-$quasiconvex function $G:[\a,\b]\rightarrow\R$ 
is integrable if $\eta$ is bounded from above on 
$G([\a,\b])\times G([\a,\b])$ (see \cite[Remark~4]{DS}).
By taking $\eta(x,y)=x-y$ in Definition~\ref{eta}, 
one recovers the classical definition of quasiconvexity. 
It is also important to note that any convex function is $\eta$-quasiconvex 
with respect to $\eta(x,y)=x-y$.  For some results around 
this recent class of functions, we invite the interested 
reader to see \cite{KKA,GDD,AwanR} and references therein.

Motivated by the above results, it is our purpose
to generalize Theorems~\ref{thm:01} and \ref{thm:Liu} 
for the class of $\eta-$quasiconvex functions. 
To the best of our knowledge, the results we prove here 
(see Theorems~\ref{MR1} and \ref{MR2}) are novel and provide
an interesting contribution to the literature of Ostrowski type results. 
In addition, we apply our results to some special known 
means of positive real numbers.

The paper is organized as follows. We begin by
recalling in Section~\ref{sec:preliminar} two  
results, needed in the sequel. In Section~\ref{MR}, we formulate 
and prove our main results, that is, Theorems~\ref{MR1} and \ref{MR2}, 
followed by several interesting corollaries. Section~\ref{app} contains 
applications of our results to special means, in particular to
the arithmetic, geometric, harmonic, logarithmic, the generalized log-mean,
and identric means (see Propositions~\ref{prop:01}, 
\ref{prop:02} and \ref{prop:03}). We end with Section~\ref{sec:conc}
of conclusion.


\section{Preliminaries}
\label{sec:preliminar}

In this section, we recall two results that will be needed 
in the proof of our main results.

\begin{lem}[See \cite{Liu}]
\label{lem1}
Let $I\subset\R$ and $G:I\rightarrow\R$ be a twice differentiable function 
on $I^{\circ}$ such that $G''\in L^1[\a,\b]$, where $\alpha,\beta\in I$ 
with $\a<\b$. Then, 
\begin{align*}
\frac{1}{\b-\a}\int^{\b}_{\a}G(x)\,dx
&-(1-\l)G\left(\frac{\a+\b}{2}\right)-\l\frac{G(\a)+G(\b)}{2}\\
&=\frac{(\b-\a)^2}{16}\Bigg[\int_{0}^{1}(x^2-2\l x)G''\left(x\frac{\a+\b}{2}+(1-x)\a\right)\,dx\\
&\qquad\qquad + \int_0^1 (x^2-2\l x)G''\left(x\frac{\a+\b}{2}+(1-x)\b\right)\,dx \Bigg]
\end{align*}
holds for any $\l\in[0, 1]$.
\end{lem}

\begin{lem}[See \cite{SA}]
\label{lem2}
Let $I\subset\R$ and $G:I\rightarrow\R$ be a twice differentiable function 
on $I^{\circ}$ such that $G''\in L^1[\a,\b]$, where $\alpha,\beta\in I$ 
with $\a<\b$. Then, 
\begin{multline*}
\frac{1}{\b-\a}\int^{\b}_{\a}G(x)\,dx-(1-\l)
G\left(\frac{\a+\b}{2}\right)-\l\frac{G(\a)+G(\b)}{2}\\
=(\b-\a)^2\int_{0}^{1}p(x)G''(x\a+(1-x)\b)\,dx
\end{multline*}
holds for any $\l\in[0, 1]$, where
\begin{equation}
\label{p}
p(x)=
\begin{cases}
\frac{1}{2}x(x-\l), & 0\leq x\leq \frac{1}{2},\\
\frac{1}{2}(1-x)(1-\l -x), & \frac{1}{2}\leq x\leq 1.
\end{cases}
\end{equation}
\end{lem}


\section{Main results}
\label{MR}

We now state and prove our first main result.
\begin{thm}
\label{MR1}
Let $I\subset[0, \infty)$ and $G:[\a, \b]\subset I\rightarrow\R$ be a twice 
differentiable function on $(\a, \b)$ with $\a<\b$. If $|G''|^q$, $q\geq 1$, 
is $\eta$-quasiconvex on $[\a, \b]$ and $\eta$-bounded from above on 
$|G''|^q([\a,\b])\times |G''|^q([\a,\b])$, then 
\begin{align*}
\Bigg|\frac{1}{\b-\a}
&\int^{\b}_{\a}G(x)\,dx-(1-\l)G\left(\frac{\a+\b}{2}\right)-\l\frac{G(\a)+G(\b)}{2}\Bigg|\\
&\leq 
\begin{cases}
\frac{(\b-\a)^2}{16}\Big(\frac{8\l^3-3\l+1}{3}\Big)\Big(
{\mathcal{N}}_{q,\eta}^{\frac{1}{q}}+{\mathcal{M}}_{q,\eta}^{\frac{1}{q}}\Big),
\quad {\rm if}~~0\leq \l\leq\frac{1}{2},\\
\frac{(\b-\a)^2}{16}\Big(\l-\frac{1}{3}\Big)\Big({\mathcal{N}}_{q,\eta}^{\frac{1}{q}}
+{\mathcal{M}}_{q,\eta}^{\frac{1}{q}}\Big), 
\quad {\rm if}~~ \frac{1}{2}\leq \l\leq 1,
\end{cases}
\end{align*}
holds, where
\begin{equation*}
{\mathcal{M}}_{q,\eta}:=\max\left\{\left|G''(\a)\right|^q,\left|G''(\a)\right|^q
+\eta\left(\left|G''\left(\frac{\a+\b}{2}\right)\right|^q, 
\left|G''(\a)\right|^q\right)\right\}
\end{equation*}
and 
\begin{equation*}
{\mathcal{N}}_{q,\eta}:=\max\left\{\left|G''(\b)\right|^q,\left|G''(\b)\right|^q
+\eta\left(\left|G''\left(\frac{\a+\b}{2}\right)\right|^q, 
\left|G''(\b)\right|^q\right)\right\}.
\end{equation*}
\end{thm}

\begin{proof}
The hypothesis that function $|G''|^q$, $q\geq 1$, is $\eta$-quasiconvex 
on $[\a, \b]$, implies that for $x\in [0,1]$ we have
\begin{equation}
\label{eta1}
\left|G''\left(x\frac{\a+\b}{2}+(1-x)\a\right)\right|^q
\leq {\mathcal{M}}_{q,\eta}
\end{equation}
and 
\begin{equation}
\label{eta2}
\left|G''\left(x\frac{\a+\b}{2}+(1-x)\b\right)\right|^q
\leq {\mathcal{N}}_{q,\eta}.
\end{equation}
Using Lemma~\ref{lem1}, \eqref{eta1}, \eqref{eta2} 
and H\"older's inequality, one obtains that
\begin{align}
\label{MI1}
\Bigg|&\frac{1}{\b-\a}\int^{\b}_{\a}G(x)\,dx-(1-\l)
G\left(\frac{\a+\b}{2}\right)-\l\frac{G(\a)+G(\b)}{2}\Bigg|\nonumber\\
&\leq \frac{(\b-\a)^2}{16}\Bigg[\int_{0}^{1}|x^2-2\l x|\Bigg|
G''\left(x\frac{\a+\b}{2}+(1-x)\a\right)\Bigg|\,dx\nonumber\\
&~~~~~~~~~~~~~~~~~~+ \int_0^1 |x^2-2\l x|\Bigg|G''\left(x\frac{\a+\b}{2}
+(1-x)\b\right)\Bigg|\,dx \Bigg]\nonumber\\
&\leq \frac{(\b-\a)^2}{16}\Bigg[\Bigg(\int_{0}^{1}|x^2-2\l x|\,dx\Bigg)^{1
-\frac{1}{q}}\nonumber\\
&\quad \times \Bigg(\int_0^1|x^2-2\l x|\Big|G''\left(x\frac{\a+\b}{2}
+(1-x)\a\right)\Big|^{q}\,dx\Bigg)^{\frac{1}{q}}\nonumber\\
&\quad + \Bigg(\int_{0}^{1}|x^2-2\l x|\,dx\Bigg)^{1-\frac{1}{q}}\Bigg(\int_0^1|x^2-2\l
x|\Big|G''\left(x\frac{\a+\b}{2}+(1-x)\b\right)\Big|^{q}\,dx\Bigg)^{\frac{1}{q}} \Bigg]\nonumber\\
&\leq \frac{(\b-\a)^2}{16}\Bigg[\Bigg(\int_{0}^{1}|x^2-2\l x|\,dx\Bigg)^{1
-\frac{1}{q}}\Bigg(\int_0^1|x^2-2\l x|{\mathcal{M}}_{q,\eta}\,dx\Bigg)^{\frac{1}{q}}\nonumber\\
&\quad+ \Bigg(\int_{0}^{1}|x^2-2\l x|\,dx\Bigg)^{1-\frac{1}{q}}\Bigg(\int_0^1|x^2-2\l 
x|{\mathcal{N}}_{q,\eta}\,dx\Bigg)^{\frac{1}{q}} \Bigg]\nonumber\\
&\leq \frac{(\b-\a)^2}{16}\Bigg[\Big({\mathcal{N}}_{q,\eta}^{\frac{1}{q}}
+{\mathcal{M}}_{q,\eta}^{\frac{1}{q}}\Big)\int_{0}^{1}|x^2-2\l x|\,dx\Bigg].
\end{align}
To finish the proof, we need to evaluate 
$\displaystyle \int_{0}^{1}|x^2-2\l x|\,dx.$ For this, we consider two cases.

{\bf Case I:} $0\leq \l\leq\frac{1}{2}$. We get $0\leq 2\l\leq 1$ and
\begin{equation}
\label{1}
\begin{aligned}
\int_{0}^{1}|x^2-2\l x|\,dx
&=\int_{0}^{2\l}|x^2-2\l x|\,dx+ \int_{2\l}^{1}|x^2-2\l x|\,dx\\
&=\int_{0}^{2\l}(2\l x-x^2)\,dx+ \int_{2\l}^{1}(x^2-2\l x)\,dx\\
&=\frac{8\l^3-3\l+1}{3}.
\end{aligned} 
\end{equation}

{\bf Case II:} $\frac{1}{2}\leq \l\leq 1$. We get $2\l\geq 1$ and 
$x^2\leq 2\l x^2\leq 2\l x$, because $x\in[0, 1]$. It follows that 
\begin{equation}
\label{2}
\begin{aligned}
\int_{0}^{1}|x^2-2\l x|\,dx
&=\int_{0}^{1}(2\l x-x^2)\,dx\\
&=\l-\frac{1}{3}.
\end{aligned} 
\end{equation}
The desired inequalities are obtained by using \eqref{1} 
and \eqref{2} in inequality \eqref{MI1}.
\end{proof}

\begin{cor}
\label{c1}
Let $I\subset[0, \infty)$ and $G:[\a, \b]\subset I\rightarrow\R$ 
be a twice differentiable function on $(\a, \b)$ with $\a<\b$. 
If $|G''|$ is $\eta$-quasiconvex on $[\a, \b]$ and $\eta$ bounded 
from above on $|G''|([\a,\b])\times |G''|([\a,\b])$, then the inequality 
\begin{align*}
\Bigg|&\frac{1}{\b-\a}\int^{\b}_{\a}G(x)\,dx
-(1-\l)G\left(\frac{\a+\b}{2}\right)-\l\frac{G(\a)+G(\b)}{2}\Bigg|\\
&\qquad \leq 
\begin{cases}
\frac{(\b-\a)^2}{16}\Big(\frac{8\l^3-3\l+1}{3}\Big)\Big({\mathcal{N}}_{\eta}
+{\mathcal{M}}_{\eta}\Big), & {\rm if}~~0\leq \l\leq\frac{1}{2},\\
\frac{(\b-\a)^2}{16}\Big(\l-\frac{1}{3}\Big)\Big({\mathcal{N}}_{\eta}
+{\mathcal{M}}_{\eta}\Big), & {\rm if}~~ \frac{1}{2}\leq \l\leq 1,
\end{cases}
\end{align*}
holds, where
\begin{equation*}
{\mathcal{M}}_{\eta}:=\max\left\{\left|G''(\a)\right|,\left|G''(\a)\right|
+\eta\left(\left|G''\left(\frac{\a+\b}{2}\right)\right|, 
\left|G''(\a)\right|\right)\right\},
\end{equation*}
and 
\begin{equation*}
{\mathcal{N}}_{\eta}:=\max\left\{\left|G''(\b)\right|,\left|G''(\b)\right|
+\eta\left(\left|G''\left(\frac{\a+\b}{2}\right)\right|, 
\left|G''(\b)\right|\right)\right\}.
\end{equation*}
\end{cor}

\begin{proof}
The proof follows by setting $q=1$ in Theorem~\ref{MR1}.
\end{proof}

\begin{rem}
By choosing different values of $\l \in [0, 1]$ in the inequality 
of Corollary~\ref{c1}, we obtain different results 
for $\eta-$quasiconvex functions. For example,
\begin{enumerate}

\item[\rm 1.] for $\l=0,$ we get a midpoint type inequality:
\begin{equation}
\label{mid1}
\Bigg|\frac{1}{\b-\a}\int^{\b}_{\a}G(x)\,dx
-G\left(\frac{\a+\b}{2}\right)\Bigg|\leq \frac{(\b-\a)^2}{48}\Big(
{\mathcal{N}}_{\eta}+{\mathcal{M}}_{\eta}\Big);
\end{equation}

\item[\rm 2.] for $\l=\frac{1}{3}$, 
we get a Simpson type inequality:
\begin{equation}
\label{sim1}
\Bigg|\frac{1}{\b-\a}\int^{\b}_{\a}G(x)\,dx
-\frac{2}{3}G\left(\frac{\a+\b}{2}\right)-\frac{G(\a)+G(\b)}{6}\Bigg|\\
\leq\frac{(\b-\a)^2}{162}\Big({\mathcal{N}}_{\eta}+{\mathcal{M}}_{\eta}\Big);
\end{equation}

\item[\rm 3.] for $\l=\frac{1}{2}$, 
we obtain a midpoint-trapezoid type inequality:
\begin{equation}
\label{midtrap1}
\Bigg|\frac{1}{\b-\a} \int^{\b}_{\a}G(x)\,dx
-\frac{1}{2}G\left(\frac{\a+\b}{2}\right)-\frac{G(\a)+G(\b)}{4}\Bigg|\\
\leq \frac{(\b-\a)^2}{96}\Big({\mathcal{N}}_{\eta}+{\mathcal{M}}_{\eta}\Big);
\end{equation}

\item[\rm 4.] for $\l=1$, we have a trapezoid type inequality:
\begin{equation}
\label{trap1}
\Bigg|\frac{1}{\b-\a}\int^{\b}_{\a}G(x)\,dx-\frac{G(\a)+G(\b)}{2}\Bigg|
\leq\frac{(\b-\a)^2}{24}\Big({\mathcal{N}}_{\eta}+{\mathcal{M}}_{\eta}\Big).
\end{equation}
\end{enumerate}
\end{rem}

Follows the second main result of our paper.

\begin{thm}
\label{MR2}
Let $I\subset[0, \infty)$ and $G:[\a, \b]\subset I\rightarrow\R$ 
be a twice differentiable function on $(\a, \b)$ with $\a<\b$. 
If $|G''|^q$, $q\geq 1$, is $\eta$-quasiconvex on $[\a, \b]$ 
and $\eta$ bounded from above on $|G''|^q([\a,\b])\times |G''|^q([\a,\b])$, 
then the inequality 
\begin{align*}
\Bigg|&\frac{1}{\b-\a}\int^{\b}_{\a}G(x)\,dx
-(1-\l)G\left(\frac{\a+\b}{2}\right)-\l\frac{G(\a)+G(\b)}{2}\Bigg|\\
&\qquad \leq 
\begin{cases}
\frac{(\b-\a)^2}{8}\Big(\frac{8\l^3-3\l+1}{3}\Big){\mathcal{U}}_{q,\eta}^{\frac{1}{q}},
& {\rm if}~~0\leq \l\leq\frac{1}{2},\\
\frac{(\b-\a)^2}{8}\Big(\frac{3\l-1}{3}\Big){\mathcal{U}}_{q,\eta}^{\frac{1}{q}},
& {\rm if}~~ \frac{1}{2}\leq \l\leq 1,
\end{cases}
\end{align*}
holds, where
\begin{equation*}
{\mathcal{U}}_{q,\eta}:=\max\left\{\left|G''(\b)\right|^q,
\left|G''(\b)\right|^q+\eta\left(\left|G''\left(\a\right)\right|^q, 
\left|G''(\b)\right|^q\right)\right\}.
\end{equation*}
\end{thm}

\begin{proof}
Since $|G''|^q$ is $\eta$-quasiconvex on $[\a, \b]$, the inequality 
\begin{equation}
\label{eta12}
\left|G''\left(x\a+(1-x)\b\right)\right|^q\leq {\mathcal{U}}_{q,\eta}
\end{equation}
holds for $x\in[0, 1]$. From the definition of $p(x)$ given by \eqref{p}, 
we observe that for $0\leq \l \leq \frac{1}{2}$ one has
\begin{equation}
\label{3}
\begin{aligned}
\int_0^1|p(x)|\,dx
&=\int_0^{\frac{1}{2}}\Big|\frac{1}{2}x(x-\l)\Big|\,dx
+\int_{\frac{1}{2}}^1\Big|\frac{1}{2}(1-x)(1-\l-x)\Big|\,dx\\
&=\frac{1}{2}\Bigg[\int_0^{\l}x(\l-x)\,dx+\int_{\l}^{\frac{1}{2}}x(x-\l)\,dx
+\int_{\frac{1}{2}}^{1-\l}(1-x)(1-\l-x)\,dx\\
&\qquad\qquad +\int_{1-\l}^1(1-x)(x-1+\l)\,dx\Bigg]\\
&=\frac{8\l^3-3\l+1}{24}.
\end{aligned}
\end{equation}
Also, for $\frac{1}{2}\leq \l\leq 1$, we get
\begin{equation}
\label{4}
\begin{aligned}
\int_0^1|p(x)|\,dx
&=\int_0^{\frac{1}{2}}\Big|\frac{1}{2}x(x-\l)\Big|\,dx
+\int_{\frac{1}{2}}^1\Big|\frac{1}{2}(1-x)(1-\l-x)\Big|\,dx\\
&=\frac{1}{2}\Bigg[\int_{0}^{\frac{1}{2}}x(\l-x)\,dx
+\int_{\frac{1}{2}}^{1}(1-x)(\l+x-1)\,dx\Bigg]\\
&=\frac{3\l-1}{24}.
\end{aligned}
\end{equation}
Now, using Lemma~\ref{lem2}, the H\"older inequality 
and \eqref{eta12}, we obtain that
\begin{align*}
&\Bigg|\frac{1}{\b-\a}\int^{\b}_{\a}G(x)\,dx-(1-\l)
G\left(\frac{\a+\b}{2}\right)-\l\frac{G(\a)+G(\b)}{2}\Bigg|\\
&\leq(\b-\a)^2\int_{0}^{1}|p(x)|\Big|G''(x\a+(1-x)\b)\Big|\,dx\\
&\leq (\b-\a)^2\Bigg(\int_{0}^{1}|p(x)|\,dx\Bigg)^{1-\frac{1}{q}}\Bigg(
\int_0^1|p(x)|\Big|G''(x\a+(1-x)\b)\Big|^q\,dx\Bigg)^{\frac{1}{q}}\\
&\leq (\b-\a)^2{\mathcal{U}}_{q,\eta}^{\frac{1}{q}}\int_{0}^{1}|p(x)|\,dx.
\end{align*}
We get the intended result by using \eqref{3} 
and \eqref{4}.
\end{proof}

\begin{cor}
\label{c2}
Let $I\subset[0, \infty)$ and $G:[\a, \b]\subset I\rightarrow\R$ 
be a twice differentiable function on $(\a, \b)$ with $\a<\b$. 
If $|G''|$ is $\eta$-quasiconvex on $[\a, \b]$ and $\eta$ bounded 
from above on $|G''|([\a,\b])\times |G''|([\a,\b])$, then the inequality 
\begin{align*}
\Bigg|&\frac{1}{\b-\a}\int^{\b}_{\a}G(x)\,dx-(1-\l)
G\left(\frac{\a+\b}{2}\right)-\l\frac{G(\a)+G(\b)}{2}\Bigg|\\
&\qquad \leq 
\begin{cases}
\frac{(\b-\a)^2}{8}\Big(\frac{8\l^3-3\l+1}{3}\Big){\mathcal{U}}_{\eta},
& {\rm if}~~0\leq \l\leq\frac{1}{2},\\
\frac{(\b-\a)^2}{8}\Big(\frac{3\l-1}{3}\Big){\mathcal{U}}_{\eta},
& {\rm if}~~ \frac{1}{2}\leq \l\leq 1,
\end{cases}
\end{align*}
holds, where
\begin{equation*}
{\mathcal{U}}_{\eta}:=\max\left\{\left|G''(\b)\right|,
\left|G''(\b)\right|+\eta\left(\left|G''\left(\a\right)\right|, 
\left|G''(\b)\right|\right)\right\}.
\end{equation*}
\end{cor}

\begin{proof}
Let $q=1$ in Theorem~\ref{MR2}.
\end{proof}

\begin{rem}
Choosing different values of $\l\in[0,1]$, we obtain, 
from Corollary~\ref{c2}, the succeeding results:
\begin{enumerate}
\item[\rm 1.] for $\l=0$, we get
\begin{align}
\label{mid2}
\Bigg|\frac{1}{\b-\a}\int^{\b}_{\a}G(x)\,dx
-G\left(\frac{\a+\b}{2}\right)\Bigg|
\leq\frac{(\b-\a)^2}{24}{\mathcal{U}}_{\eta};
\end{align}

\item[\rm 2.] for $\l=\frac{1}{3}$, we obtain
\begin{align}
\Bigg|\frac{1}{\b-\a}\int^{\b}_{\a}G(x)\,dx
-\frac{2}{3}G\left(\frac{\a+\b}{2}\right)-\frac{G(\a)+G(\b)}{6}\Bigg|
\leq \frac{(\b-\a)^2}{81}{\mathcal{U}}_{\eta};
\end{align}

\item[\rm 3.] for $\l=\frac{1}{2}$, we have
\begin{align}
\Bigg|\frac{1}{\b-\a}\int^{\b}_{\a}G(x)\,dx
-\frac{1}{2}G\left(\frac{\a+\b}{2}\right)-\frac{G(\a)
+G(\b)}{4}\Bigg|\leq \frac{(\b-\a)^2}{48}{\mathcal{U}}_{\eta};
\end{align}

\item[\rm 4.] for $\l=1$, we get
\begin{align}\label{trap2}
\Bigg|\frac{1}{\b-\a}\int^{\b}_{\a}G(x)\,dx-\frac{G(\a)+G(\b)}{2}\Bigg|
\leq \frac{(\b-\a)^2}{12}{\mathcal{U}}_{\eta}.
\end{align}
\end{enumerate}
\end{rem}

\begin{rem}
Let $0<\a<\b$. By setting $G(x)=\ln x$ with $x\in[\a,\b]$ and $\eta(x,y)=x-y$ 
in inequalities \eqref{mid2}--\eqref{trap2}, one gets {\rm \cite[Proposition 3]{Liu}}.
\end{rem}


\section{Application to special means}
\label{app}

In this section, we apply our results to the following special means 
of arbitrary positive numbers $\mu$ and $\nu$ with $\mu\neq\nu$:
\begin{enumerate}
\item the arithmetic mean
$$
A(\mu,\nu)=\frac{\mu+\nu}{2};
$$

\item the geometric mean
$$
G(\mu,\nu)=\sqrt{\mu\nu};
$$
\item the harmonic mean 
$$
H(\mu,\nu)=\frac{2\mu\nu}{\mu+\nu};
$$

\item the logarithmic mean
$$
L(\mu,\nu)=\frac{\nu-\mu}{\ln\nu-\ln\mu};
$$

\item the generalized log-mean
$$
L_p(\mu,\nu)=\Bigg[\frac{\nu^{p+1}-\mu^{p+1}}{(p+1)(\nu
-\mu)}\Bigg]^{\frac{1}{p}}, \quad p\neq -1, 0;
$$

\item the identric mean
$$
I(\mu,\nu)=\frac{1}{e}\Bigg(\frac{\nu^{\nu}}{\mu^{\mu}}\Bigg)\frac{1}{\nu-\mu}.
$$
\end{enumerate}

We now state our findings in the following propositions.

\begin{pro}
\label{prop:01}
Let $\mu$ and $\nu$ be two positive numbers, 
$\mu<\nu$. The following inequalities hold:
\begin{enumerate}
\item $\left|L^2_2(\mu,\nu)-A^2(\mu,\nu)\right|
\leq\displaystyle \frac{(\nu-\mu)^2}{12}$;

\item $\left|L^2_2(\mu,\nu)-\frac{2A^2(\mu,\nu)+A(\mu^2,\nu^2)}{3}\right|
\leq\displaystyle\frac{2(\nu-\mu)^2}{81}$;

\item $\left|L^2_2(\mu,\nu)-\frac{A^2(\mu,\nu)
+A(\mu^2,\nu^2)}{2}\right|
\leq\displaystyle \frac{(\nu-\mu)^2}{24}$;

\item $\left|L^2_2(\mu,\nu)-A(\mu^2,\nu^2)\right|
\leq\displaystyle\frac{(\nu-\mu)^2}{6}$.
\end{enumerate}
\end{pro}

\begin{proof}
The desired inequalities follow by employing \eqref{mid1}--\eqref{trap1} 
to function $G(x)=x^2$ defined on the interval $[\mu,\nu]$. 
In this case, $|G''(x)|=2$. By taking $\eta(x,y)=x-y$, we easily see that 
$|G''(x)|$ is $\eta-$quasiconvex. Moreover, 
${\mathcal{M}}_{\eta}={\mathcal{N}}_{\eta}=2$.
\end{proof}

\begin{pro}
\label{prop:02}
Let $\mu$ and $\nu$ be two positive numbers, $\mu<\nu$. 
The following inequalities hold:
\begin{enumerate}
\item $\left|A^{-1}(\mu,\nu)-L^{-1}(\mu,\nu)\right|
\leq\frac{(\nu-\mu)^2}{48}\Bigg[\max\left\{\frac{2}{\mu^3}, 
\frac{16}{(\mu+\nu)^3}\right\}+\max\left\{\frac{2}{\nu^3}, 
\frac{16}{(\mu+\nu)^3}\right\}\Bigg]$;

\item $\left|\frac{2A^{-1}(\mu,\nu)+H^{-1}(\mu,\nu)}{3}
-L^{-1}(\mu,\nu)\right|\leq \frac{(\nu-\mu)^2}{162}\Bigg[
\max\left\{\frac{2}{\mu^3}, \frac{16}{(\mu+\nu)^3}\right\}
+\max\left\{\frac{2}{\nu^3}, \frac{16}{(\mu+\nu)^3}\right\}\Bigg]$;

\item $\left|\frac{A^{-1}(\mu,\nu)+H^{-1}(\mu,\nu)}{2}
-L^{-1}(\mu,\nu)\right|\leq \frac{(\nu-\mu)^2}{96}\Bigg[
\max\left\{\frac{2}{\mu^3}, \frac{16}{(\mu+\nu)^3}\right\}
+\max\left\{\frac{2}{\nu^3}, \frac{16}{(\mu+\nu)^3}\right\}\Bigg]$;

\item $\left|H^{-1}(\mu,\nu)-L^{-1}(\mu,\nu)\right|
\leq\frac{(\nu-\mu)^2}{24}\Bigg[\max\left\{\frac{2}{\mu^3}, 
\frac{16}{(\mu+\nu)^3}\right\}+\max\left\{\frac{2}{\nu^3}, 
\frac{16}{(\mu+\nu)^3}\right\}\Bigg]$.
\end{enumerate}
\end{pro}

\begin{proof}
We apply inequalities \eqref{mid1}--\eqref{trap1} to the function 
$G:[\mu,\nu]\rightarrow\mathbb{R}$ defined by $G(x)=\frac{1}{x}$. 
For this, we observe that $|G''(x)|=\frac{2}{x^3}$ is convex on 
$[\mu,\nu]$ and so $\eta-$quasiconvex with respect 
to $\eta(x,y)=x-y$.
\end{proof}

We end with more four new inequalities.

\begin{pro}
\label{prop:03}
Let $\mu$ and $\nu$ be two positive numbers with $\mu<\nu$. 
Then the following inequalities hold:
\begin{enumerate}
\item $\left|\ln A(\mu,\nu)-\ln I(\mu,\nu)\right|
\leq\frac{(\nu-\mu)^2}{48}\Bigg[\max\left\{\frac{1}{\mu^2}, 
\frac{4}{(\mu+\nu)^2}\right\}+\max\left\{\frac{1}{\nu^2}, 
\frac{4}{(\mu+\nu)^2}\right\}\Bigg]$;

\item $\left|\frac{2\ln A(\mu,\nu)+\ln G(\mu,\nu)}{3}
-\ln I(\mu,\nu)\right|\leq \frac{(\nu-\mu)^2}{162}\Bigg[
\max\left\{\frac{1}{\mu^2}, \frac{4}{(\mu+\nu)^2}\right\}
+\max\left\{\frac{1}{\nu^2}, \frac{4}{(\mu+\nu)^2}\right\}\Bigg]$;

\item $\left|\frac{\ln A(\mu,\nu)+\ln G(\mu,\nu)}{2}
-\ln I(\mu,\nu)\right|\leq \frac{(\nu-\mu)^2}{96}\Bigg[
\max\left\{\frac{1}{\mu^2}, \frac{4}{(\mu+\nu)^2}\right\}
+\max\left\{\frac{1}{\nu^2}, \frac{4}{(\mu+\nu)^2}\right\}\Bigg]$;

\item $\left|\ln G(\mu,\nu)-\ln I(\mu,\nu)\right|
\leq\frac{(\nu-\mu)^2}{24}\Bigg[\max\left\{\frac{1}{\mu^2}, 
\frac{4}{(\mu+\nu)^2}\right\}+\max\left\{\frac{1}{\nu^2}, 
\frac{4}{(\mu+\nu)^2}\right\}\Bigg]$.
\end{enumerate}
\end{pro}

\begin{proof}
Result follows by applying \eqref{mid1}--\eqref{trap1} 
to the function $G(x)=\ln x$, $x\in [\mu,\nu]$, 
taking $\eta(x,y)=x-y$ and noting that 
$|G''(x)|=\frac{1}{x^2}$ is $\eta-$quasiconvex.
\end{proof}


\section{Conclusion}
\label{sec:conc}

We proved two main theorems that establish Ostrowski type
inequalities in terms of a parameter $\l\in [0, 1]$. By choosing 
$\l=0$, $1/3$, $1/2$, $1$, we deduced midpoint, Simpson, 
midpoint-trapezoid and trapezoid type inequalities, respectively.
Thereafter, we illustrated the importance of our results 
by applying them to special means of positive real numbers.


\begin{acknowledgement}
This research was partially supported by the Portuguese Foundation
for Science and Technology (FCT) through CIDMA, 
project UID/MAT/04106/2013.
\end{acknowledgement}



\end{document}